\documentclass[twoside]{article}
\usepackage{mathrsfs}
\usepackage{amsfonts}
\usepackage{stmaryrd}
\usepackage{graphicx,amssymb,mathrsfs,amsmath,latexsym,amsfonts,amsthm,bbding}

\renewcommand{\baselinestretch}{1.5} 
\catcode`@=11 \long\def\@makefntext#1{\noindent #1}
\newskip\tabcentering \tabcentering=1000pt plus 1000pt minus 1000pt

\def\MCH#1#2{\setbox0=\hbox{\raise#1\hbox{#2}}\smash{\box0}}

\let\@oddfoot\@empty  \let\@evenfoot\@empty

\def\@evenfoot{\footnotesize{\thepage\hfill}}\def\@oddfoot{\footnotesize{\hfill\thepage}}
\def\@evenhead{\hfill\footnotesize{{\it Zhao Shien} \&
{\it Shi guang}}}
\def\@oddhead{{\footnotesize
{\it A geometric form of Hahn-Banach extension theorem for $L^{0}-$
linear functions } \hfill}}

\def\leq{\leqslant}
\def\geq{\geqslant}
       
  \def\wj{\end{document}}  \newsymbol\wjzhml 203F \def\no{\noindent}

\begin{document}
\abovedisplayskip=3pt plus 1pt minus 2pt 
\belowdisplayskip=3pt plus 1pt minus 2pt 

\

\textwidth=145truemm \textheight=218truemm \thispagestyle{empty}
\renewcommand{\baselinestretch}{1.5}\baselineskip 9pt

\vspace{8true mm}

\renewcommand{\baselinestretch}{1.5}\baselineskip 15pt

\noindent{\Large\bf A geometric form of the Hahn-Banach extension
theorem for $L^{0}-$ linear functions and the Goldstine-Weston
theorem in random normed modules}퉀
\vspace{0.5 true cm}

\noindent{\normalsize\sf Zhao ShiEn$^{\ast}$, Shi Guang 
\footnotetext{\baselineskip 10pt
$^{\textmd{\FiveStarOpen}}$Supported by the National Natural
Science Foundation of China (No. 10871016).\\%
$^{\ast}$ Corresponding author. \\
E-mail addresses: zsefh@ss.buaa.edu.cn, g\_shi@ss.buaa.edu.cn.}}

\vspace{0.2 true cm} \baselineskip 12pt \noindent{\footnotesize\rm
LMIB and School of Mathematics and Systems Science, Beihang
University, Beijing 100191, PR China \\

\noindent{{\bf\small Abstract}\hspace{2.8mm} ~~{\small  In this
paper, we present a geometric form of the Hahn-Banach extension
theorem for $L^{0}-$linear functions and prove that the geometric
form is equivalent to the analytic form of the Hahn-Banach extension
theorem. Further, we use the geometric form to give a new proof of a known
basic strict separation theorem in random locally convex modules.
Finally, using the basic strict separation theorem we establish the
Goldstine-Weston theorem in random normed modules under the two kinds of
topologies----the $(\varepsilon,\lambda)-$topology
and the locally $L^{0}-$convex topology, and also provide a counterexample showing that the Goldstine-Weston theorem
under the locally $L^{0}-$convex topology can only hold for random normed modules with the countable concatenation property.
    } }

\vspace{1mm} \no{\footnotesize{\bf Keywords:\hspace{2mm} Hahn-Banach
extension theorem, random locally convex module, random normed
module, $(\epsilon,\lambda)-$topology, locally $L^{0}-$convex
topology, separation theorem, Goldstine-Weston theorem

\no{\footnotesize{\bf MSC(2000):\hspace{2mm} 46A22, 46A16, 46H25, 46H05}
 \vspace{2mm}

\renewcommand{\baselinestretch}{1.48}
\parindent=16pt  \parskip=2mm
\rm\normalsize\rm

\section{Introduction}

It is well known that the classical Hahn-Banach extension theorem for linear functionals has both its algebraic form and geometric form. The corresponding algebraic form of the Hahn-Banach extension
theorem for random linear functionals are due to Guo in
\cite{Guo-Survey I,Guo-J}. The Hahn-Banach extension
theorem for
$L^{0}-$linear functions, namely Proposition 1.1 below, is due to \cite{B-S,{VD}}, an extremely simple proof of which was given in
\cite{Guo-Relations between}.

Before giving Proposition 1.1, we first recall some notation and
terminology.

In the sequel of this paper, $(\Omega,\mathcal {F},P)$ denotes a
probability space, $N$ the set of all positive integers, $K$ the
real number field $R$ or the complex number field $C$,
$\bar{R}=[-\infty,+\infty]$, $\bar{L}^{0}(\mathcal {F},R)$ the set
of equivalence classes of extended real-valued random variables on $(\Omega,\mathcal {F},P)$,
$L^{0}(\mathcal {F},K)$ the algebra of equivalence classes of
$K-$valued random variables on $(\Omega,\mathcal {F},P)$
under the ordinary scalar multiplication, addition and
multiplication operations on equivalence classes, the null and unit
elements are still denoted by $0$ and $1$, respectively.

It is well known from \cite{DS} that $\bar{L}^{0}(\mathcal {F},R)$
is a complete lattice under the ordering $\leq$: $\xi\leq\eta$ iff
$\xi^{0}(\omega)\leq\eta^{0}(\omega)$, for almost all $\omega$ in
$\Omega$ (briefly, a.s.), where $\xi^{0}$ and $\eta^{0}$ are
arbitrarily chosen representatives of $\xi$ and $\eta$, respectively
(see also Proposition 2.1 below). Furthermore, every subset
$G$ of $\bar{L}^{0}(\mathcal {F},R)$ has a supremum and an infimum, denoted by
$\bigvee G$ and $\bigwedge G$, respectively. In particular,
$L^{0}(\mathcal {F},R)$, as a sublattice of $\bar{L}^{0}(\mathcal
{F},R)$, is also a complete lattice in the sense that every subset
with an upper bound has a supremum.

Specially, $L^{0}_{+}=\{\xi\in L^{0}(\mathcal{F},R)~|~\xi\geq 0\}$,
$L^{0}_{++}=\{\xi\in L^{0}(\mathcal{F},R)~|~\xi>0$ on $\Omega\}$,
where for $A\in \mathcal{F}$, $``\xi>\eta"$ on $A$ means
$\xi^{0}(\omega)>\eta^{0}(\omega)$ a.s. on $A$ for any chosen
representatives $\xi^{0}$ and $\eta^{0}$ of $\xi$ and $\eta$,
respectively. As usual, $\xi>\eta$ means $\xi\geq\eta$ and
$\xi\neq\eta$.

Given a random locally convex module $(E,\mathcal{P})$ over $K$ with
base $(\Omega,\mathcal {F},P)$, let
$\mathcal{T}_{\varepsilon,\lambda}$ and $\mathcal{T}_{c}$ denote the
$(\varepsilon,\lambda)-$topology and the locally $L^{0}-$convex
topology for $E$, respectively, see \cite{Guo-Relations between,FKV}
and also Section 2 for the definitions of these two kinds of
topologies.

\no{\bf Proposition 1.1 (The algebraic form of Hahn-Banach theorem
for $L^{0}-$linear functions \cite{B-S,VD,FKV}).}\quad Let $E$ be a left module over the algebra
$L^{0}(\mathcal{F},R)$, $M$ an $L^{0}(\mathcal{F},R)-$submodule in
$E$, $g:M\rightarrow L^{0}({\cal F},R)$ an $L^{0}-$linear functional
and $p:E\rightarrow L^{0}({\cal F},R)$ an $L^{0}-$sublinear
functional such that $g(x)\leq p(x),\forall x\in M$. Then there
exists an $L^{0}-$linear functional $f:E\rightarrow L^{0}({\cal
F},R)$ such that $f$ extends $g$ and $f(x)\leq p(x),\forall x\in E$.

In this paper we present the following geometric form of
Proposition 1.1, namely Proposition 1.2 below, and point out that the geometric form is equivalent
to the algebraic form stated above.

\no{\bf Proposition 1.2(The geometric form of Hahn-Banach theorem for
$L^{0}-$linear functions).}\quad Let $E$ be a left module over the
algebra $L^{0}(\mathcal{F},R)$, $M$ an
$L^{0}(\mathcal{F},R)-$submodule in $E$ and $G$ an $L^{0}-$convex
and $L^{0}-$absorbent subset of $E$. If $g:M\rightarrow L^{0}({\cal
F},R)$ is an $L^{0}-$linear functional and $g(y)\leq 1$ for any
$y\in M\cap G$, then there exists an $L^{0}-$linear functional
$f:E\rightarrow L^{0}({\cal F},R)$ such that $f$ extends $g$ and
$f(x)\leq 1,\forall x\in G$.

In addition, we make use of the geometric form
to give a new proof of the following known basic strict separation theorem
in random locally convex modules:

\no {\bf Proposition 1.3 (\cite{GXC}).}\quad {\it Let $(E, \mathcal
{P})$ be a random locally convex module over $K$ with base
$(\Omega,\mathcal {F},P)$, $G$ a $\mathcal
{T}_{\varepsilon,\lambda}-$closed and $L^{0}-$convex subset of $E$,
$x_0\in E\setminus G,\ \xi_{\mathcal
{Q}}=\bigwedge\{\|x_{0}-h\|_{\mathcal {Q}}~|~h\in G\}$ for each
$\mathcal {Q}\in {\mathcal {F}(\mathcal {P})}$ and
$\xi=\bigvee\{\xi_{\mathcal {Q}}~|~\mathcal {Q}\in {\mathcal
{F}(\mathcal {P})}\}$. Then there exists
a continuous module homomorphism $f$ from $(E,\mathcal
{T}_{\epsilon,\lambda})$ to $(L^{0}(\mathcal {F},K),\mathcal
{T}_{\epsilon,\lambda})$ such that
$$(\mbox{Re}f)(x_0)>\bigvee\{(\mbox{Re}f)(y)~|~ y\in G\},$$
where $\mbox{Re}f$ denotes the real part of $f$, namely
$f(x)=(\mbox{Re}f)(x)-i(\mbox{Re}f)(ix),\forall x\in E$ and
$$(\mbox{Re}f)(x_0)>\bigvee\{(\mbox{Re}f)(y)~|~ y\in G\}~on~[\xi>0].$$}
\no In the final part of this paper, we establish the Goldstine-Weston
theorem in random normed modules under the two kinds of topologies, namely the $(\varepsilon,\lambda)-$topology and the locally $L^{0}-$convex topology, which are stated as follows:

\no {\bf Theorem 1.1.} \quad {\it Let $(E,\|\cdot\|)$ be an $RN$
module over $K$ with base $(\Omega,\mathcal{F},P)$, $J$ the
natural embedding mapping: $E\rightarrow E^{\ast\ast}$
defined by $J(x)(g)=g(x)$ for any $x\in E$ and $g\in E^{\ast}$,
$E(1)=\{x\in E~|~\|x\|\leq 1\}$ and
$\overline{J(E(1))}_{\varepsilon,\lambda}^{w^{\ast}}$ the
closure of $J(E(1))$ with respect to
$\sigma_{\varepsilon,\lambda}(E^{\ast\ast},E^{\ast})$. Then
$\overline{J(E(1))}_{\varepsilon,\lambda}^{w^{\ast}}=E^{\ast\ast}(1)$,
where $E^{\ast\ast}(1)=\{\phi\in
E^{\ast\ast}~|~\|\varphi\|^{\ast\ast}\leq1\}$.}

\no {\bf Theorem 1.2.} \quad {\it Let $(E,\|\cdot\|)$ be an $RN$
module over $K$ with base $(\Omega,\mathcal{F},P)$ such that $E$ has
the countable concatenation property, $J$ and $E(1)$ the same as in
Theorem 1.1, and $\overline{J(E(1))}_{c}^{w^{\ast}}$ the closure
of $J(E(1))$ with respect to $\sigma_{c}(E^{\ast\ast},E^{\ast})$.
Then $\overline{J(E(1))}_{c}^{w^{\ast}}=E^{\ast\ast}(1)$. }

\no Further, we give an example to show that $J(E(1))$ may not be
dense in $E^{\ast\ast}(1)$ under $\sigma_{c}(E^{\ast\ast},E^{\ast})$
if $(E,\|\cdot\|)$ has not the countable concatenation property.

The remainder of this paper is organized as follows: in Section 2 we
will recapitulate some known basic facts, in Section 3 we will prove that the geometric form of Hahn-Banach extension theorem for
$L^{0}-$linear functions is equivalent to the algebraic form and in Section 4 we will prove the
Goldstine-Weston theorem in random normed modules.

\section{Preliminaries}

\no {\bf Proposition 2.1 (\cite{DS}).}\quad {\it For every subset
$G$ of $\bar{L}^{0}(\mathcal {F},R)$ there exist countable subsets
$\{a_{n}~|~n\in N\}$ and $\{b_{n}~|~n\in N\}$ of $G$ such that
$\bigvee G=\bigvee_{n\geq1}a_{n}$ and $\bigwedge
G=\bigwedge_{n\geq1}b_{n}$. Further, if $G$ is directed $($dually
directed$)$ with respect to $\leq$, then the above $\{a_{n}~|~n\in
N\}$ $($accordingly, $\{b_{n}~|~n\in N\}$$)$ can be chosen as
nondecreasing $($correspondingly, nonincreasing$)$ with respect to
$\leq$.}

For an arbitrarily chosen representative $\xi^{0}$ of $\xi\in
L^{0}(\mathcal{F},K)$, define the two
random variables $(\xi^{0})^{-1}$ and $|\xi^{0}|$ by
$(\xi^{0})^{-1}(\omega)=1/\xi^{0}(\omega)$ if $\xi^{0}(\omega)\neq
0$, and $(\xi^{0})^{-1}(\omega)=0$ otherwise, and by
$|\xi^{0}|(\omega)=|\xi^{0}(\omega)|$, $\forall\omega\in \Omega$.
Then the equivalent class $Q(\xi)$ of $(\xi^{0})^{-1}$ is called the
generalized inverse of $\xi$ and the equivalent class $|\xi|$ of
$|\xi^{0}|$ the absolute value of $\xi$.

Besides, for any $A\in \mathcal{F}$, $A^{c}$ denotes the complement
in $\Omega$, $\tilde{A}:=\{B\in \mathcal{F}~|~P(A\Delta
B)=0\}$ the equivalence class of $A$, where $\Delta$ is the
symmetric difference operation, $I_{A}$ the characteristic function
of $A$ and $\tilde{I}_{A}$ the equivalence class
of $I_{A}$. Given two $\xi$ and $\eta$ in $L^{0}(\mathcal{F}, R)$,
and $A=\{\omega\in\Omega~|~\xi^{0}\neq\eta^{0}\}$, where $\xi^{0}$
and $\eta^{0}$ are arbitrarily chosen representatives of $\xi$ and
$\eta$ respectively, then we always write $[\xi\neq\eta]$ for the
equivalence class of $A$ and $I_{[\xi\neq\eta]}$ for
$\tilde{I}_{A}$, one can also understand the implication of such
notations as $I_{[\xi\leq\eta]}$, $I_{[\xi<\eta]}$ and
$I_{[\xi=\eta]}$.

\no  {\bf  Definition 2.1 (\cite{{Guo-Survey},{G-C}}).}\quad
(1)~~{\it Let $E$ be a linear space over $K$, then a mapping $f:
E\rightarrow L^{0}(\mathcal {F}, K)$ is called a random linear
functional on $E$ if $f$ is linear;

$(2)$  If $E$ is a linear space over $R$, then a mapping $f:
E\rightarrow L^{0}(\mathcal {F},R)$ is called a random sublinear
functional on $E$ if $f(\alpha x)=\alpha\cdot f(x)$ for any positive
real number $\alpha$ and $x\in E$, and if $f(x+y)\leq f(x)+f(y),\forall
x,y\in E$;

$(3)$  Let $E$ be a linear space over $K$, then a mapping $f:
E\rightarrow L_{+}^{0}$ is called a random seminorm on $E$ if
$f(\alpha x)=|\alpha|\cdot f(x), \forall \alpha \in K$ and $x\in E$,
and if $f(x+y)\leq f(x)+f(y), \forall x,y\in E$;

$(4)$  Let $E$ be a left module over the algebra $L^{0}(\mathcal
{F},K)$, then a mapping $f: E\rightarrow L^{0}(\mathcal {F}, K)$ is
called a $L^{0}-$linear function on $E$ if $f$ is a module
homomorphism;

$(5)$  Let $E$ be a left module over the algebra $L^{0}(\mathcal
{F},R)$, a mapping $f: E\rightarrow L^{0}(\mathcal {F}, R)$ is
called an $L^{0}$-sublinear functional on $E$ if $f$ is a random
sublinear function on $E$ such that $f(\xi\cdot x)=\xi\cdot f(x),
\forall\xi \in L_{+}^{0}$ and $x\in E$;

$(6)$  Let $E$ be a left module over the algebra $L^{0}(\mathcal
{F},K)$, then a mapping $f: E\rightarrow L_{+}^{0}$ is called an
$L^{0}$-seminorm on $E$ if $f$ is a random seminorm on $E$ such that
$f(\xi\cdot x)=|\xi|\cdot f(x), \forall\xi \in L^{0}(\mathcal
{F},K)$ and $x\in E$. }

\no  {\bf  Definition 2.2 (\cite{{Guo-Relations
between},{Guo-Survey},{Guo-progress}}).}\quad {\it An ordered pair
$(E,\mathcal {P})$ is called a random locally convex space over $K$
with base $(\Omega , \mathcal{F}, P)$ if $E$ is a linear space over
$K$ and $\mathcal{P}$ is a family of random seminorms on $E$ such
that the following axiom is satisfied: }

(1) {\it $\bigvee\{\|x\|~|~\|\cdot\|\in{\cal P}\}=0$ implies
$x=\theta$ $($the null element of $E$ $)$.\\
In addition, if $E$ is a left module over the algebra
$L^{0}(\mathcal{F},K)$ and each $\|\cdot\|$ in $\mathcal{P}$ is an
$L^{0}-$seminorm, then such a random locally convex space is called
a random locally convex module.}

\no  {\bf Remark 2.1.}\quad Let $(E,\mathcal {P})$ be a random
locally convex space (a random locally convex module) over $K$ with
base $(\Omega , \mathcal{F}, P)$. If $\mathcal{P}$ degenerates to a
singleton $\{\|\cdot\|\}$, then $(E,\|\cdot\|)$ is exactly a random
normed space (briefly, an $RN$ space) (correspondingly, a random
normed module (briefly, an $RN$ module)). Specially,
$(L^{0}(\mathcal{F},K),|\cdot|)$ is an $RN$ module.

In the sequel, for a random locally convex space $(E,{\cal P})$ with
base $(\Omega,{\cal F},P)$ and for each finite subfamily
$\mathcal{Q}$ of ${\cal P}$, $\|\cdot\|_{\mathcal{Q}}:E\rightarrow
L^{0}_{+}({\cal F})$ always denotes the random seminorm of $E$
defined by $\|x\|_{\mathcal{Q}}=\bigvee\{\|x\|~|~\|\cdot\|\in
\mathcal{Q}\},\forall x\in E$, and ${\cal F}({\cal P})$ the set of
finite subfamilies of ${\cal P}$.

For each random locally convex space $(E,\mathcal {P})$ over $K$
with base $(\Omega , \mathcal{F}, P)$, $\mathcal {P}$ can induce the following two
kinds of topologies, namely the $(\varepsilon,\lambda)-$topology and
the locally $L^{0}-$convex topology.

\no  {\bf  Definition 2.3 (\cite{{Guo-Relations
between},{Guo-Survey},{Guo-progress}}).}\quad {\it Let
$(E,\mathcal{P})$ be a random locally convex space over $K$ with
base $(\Omega,\mathcal{F},P)$. For any positive real numbers
$\varepsilon$ and $\lambda$ such that $0<\lambda<1$, and any
$\mathcal {Q}\in \mathcal{F}(\mathcal{P})$, let $N_{\theta}(\mathcal
{Q},\varepsilon,\lambda)=\{x\in
E~|~P\{\omega\in\Omega~|~\|x\|_{\mathcal
{Q}}(\omega)<\varepsilon\}>1-\lambda\}$, then $\{N_{\theta}(\mathcal
{Q},\varepsilon,\lambda)~|~\mathcal {Q}\in
\mathcal{F}(\mathcal{P}),\varepsilon>0,0<\lambda<1\}$ is easily
verified to be a local base at the null vector $\theta$ of some
Hausdorff linear topology, called the
$(\varepsilon,\lambda)-$topology for $E$ induced by $\mathcal{P}$.}

From now on, the $(\varepsilon,\lambda)-$topology for each random
locally convex space is always denoted by $\mathcal
{T}_{\varepsilon,\lambda}$ when no confusion occurs.

\no  {\bf  Definition 2.4(\cite{FKV, Guo-progress}).}\quad {\it Let
$(E,\mathcal{P})$ be a random locally convex space over $K$ with
base $(\Omega,\mathcal{F},P)$. For any $\mathcal {Q}\in
\mathcal{F}(\mathcal{P})$ and $\varepsilon\in L^{0}_{++}$, let
$N_{\theta}(\mathcal {Q},\varepsilon)=\{x\in
E~|~\|x\|_{\mathcal{Q}}\leq \varepsilon\}$. A subset $G$ of $E$ is
called $\mathcal{T}_{c}-$open if for each $x\in G$ there exists some
$N_{\theta}(\mathcal {Q},\varepsilon)$ such that
$x+N_{\theta}(\mathcal {Q},\varepsilon)\subset G$, $\mathcal
{T}_{c}$ denotes the family of $\mathcal{T}_{c}-$open subsets of
$E$. Then it is easy to see that $(E,\mathcal {T}_{c})$ is a
Hausdorff topological group with respect to the addition on $E$.
$\mathcal{T}_{c}$ is called the locally $L^{0}-$convex topology for
$E$ induced by $\mathcal{P}$.}

From now on, the locally $L^{0}-$convex topology for each random
locally convex space is always denoted by $\mathcal{T}_{c}$ when no
confusion occurs.

Now, we present the definition of random conjugate spaces of a
random locally convex space. Historically, the earliest two notions
of a random conjugate space of a random locally convex space were
introduced in \cite{Guo-Survey, Guo-Module}, respectively. As shown
in \cite{Guo-Relations between,Guo-progress}, it turned out that
they just correspond to the $(\varepsilon,\lambda)-$topology and the
locally $L^{0}-$convex topology in the context of a random locally
convex module, respectively!

\no  {\bf  Definition 2.5(\cite{Guo-Module}).}\quad {\it Let
$(E,\mathcal{P})$ be a random locally convex space over $K$ with
base $(\Omega,\mathcal{F},P)$. A random linear functional
$f:E\rightarrow L^{0}(\mathcal{F},K)$ is called an a.s. bounded
random linear functional of type I if there are some $\xi\in
L^{0}_{+}$ and $\mathcal{Q}\in \mathcal{F}(\mathcal{P})$ such that
$|f(x)|\leq\xi\cdot\|x\|_{\mathcal{Q}},\forall x\in E$. Denote by
$E^{\ast}_{I}$ the set of a.s. bounded random linear functional of
type I on $E$. The module multiplication operation
$\cdot:L^{0}(\mathcal{F},K)\times E^{\ast}_{I}\rightarrow
E^{\ast}_{I}$ is defined by $(\xi f)(x)=\xi(f(x)),\forall \xi\in
L^{0}(\mathcal{F},K),f\in E^{\ast}_{I}$ and $x\in E$. It is easy to
see that $E^{\ast}_{I}$ is a left module over
$L^{0}(\mathcal{F},K)$, called the random conjugate space of type I
of $E$.}

\no  {\bf  Definition 2.6(\cite{G-C}).}\quad {\it Let
$(E,\mathcal{P})$ be a random locally convex space over $K$ with
base $(\Omega,\mathcal{F},P)$. A random linear functional
$f:E\rightarrow L^{0}(\mathcal{F},K)$ is called an a.s. bounded
random linear functional of type II on $E$ if there exist a
countable partition $\{A_{i}~|~i\in N\}$ of $\Omega$ to
$\mathcal{F}$, a sequence $\{\xi_{i}~|~i\in N\}$ in $L^{0}_{+}$ and
a sequence $\{\mathcal {Q}_{i}~|~i\in N\}$ in
$\mathcal{F}(\mathcal{P})$ such that
$|f(x)|\leq\Sigma_{i=1}^{\infty}\tilde{I}_{A_{i}}\cdot\xi_{i}\cdot
\|x\|_{\mathcal{Q}_{i}},\forall x\in E$. Denote by $E^{\ast}_{II}$
the $L^{0}(\mathcal{F},K)-$module of a.s. bounded random linear
functional of type II on $E$, called the random conjugate space of
type II of $E$.}

\no  {\bf  Definition 2.7.}\quad {\it Let $(E,\mathcal{P})$ be a
random locally convex module over $K$ with base
$(\Omega,\mathcal{F},P)$ and define
$E_{\varepsilon,\lambda}^{\ast}$, $E_{c}^{\ast}$ as follows:
$$(1)~~~E_{\varepsilon,\lambda}^{\ast}=\{f~|~f~is~a~continuous~module~homomorphism~from~(E,{\cal
T}_{\varepsilon,\lambda})~to~(L^{0}(\mathcal {F},K),{\cal
T}_{\varepsilon,\lambda})\},~$$
$$(2)~~~E_{c}^{\ast}=\{f~|~f~is~a~continuous~module~homomorphism~from~(E,{\cal
T}_{c})~to~(L^{0}(\mathcal {F},K),{\cal T}_{c})\}.~~~~~~$$}
~~~~Propositions 2.2 and 2.3 below give the topological
characterizations of an element in $E^{\ast}_{I}$ and
$E^{\ast}_{II}$, respectively.

\no {\bf Proposition 2.2(\cite{Guo-Relations between,
Guo-Survey}).}\quad {\it Let $(E,{\cal P})$ be a random locally
convex module over $K$ with base $(\Omega,{\cal F},P)$ and
$f:E\rightarrow L^{0}({\cal F},K)$ a random linear functional. Then
$f\in E^{\ast}_I$ iff $f$ is a continuous module homomorphism from
$(E,{\cal T}_c)$ to $(L^{0}({\cal F},K),{\cal T}_c)$, namely
$E^{\ast}_I=E^{\ast}_c$.}

\no {\bf Proposition 2.3(\cite{Guo-Relations between,
Guo-Zhu}).}\quad {\it Let $(E,{\cal P})$ be a random locally convex
module over $K$ with base $(\Omega,{\cal F},P)$ and $f:E\rightarrow
L^{0}({\cal F},K)$ a random linear functional. Then $f\in
E^{\ast}_{II}$ iff $f$ is a continuous module homomorphism from
$(E,{\cal T}_{\epsilon,\lambda})$ to $(L^{0}({\cal F},K),{\cal
T}_{\epsilon,\lambda})$, namely
$E^{\ast}_{II}=E^{\ast}_{\varepsilon,\lambda}$.}

\no {\bf Remark 2.2.}\quad It is clear that $E^{\ast}_{c}\subset
E_{\varepsilon,\lambda}^{\ast}$ from Proposition 2.2 and 2.3.
Specially, if $(E,\|\cdot\|)$ is an $RN$ module over $K$ with base
$(\Omega,\mathcal{F},P)$, then
$E^{\ast}_{c}=E_{\varepsilon,\lambda}^{\ast}$ (see
\cite{Guo-Relations between} for details), in which case we denote $E^{\ast}_{\varepsilon,\lambda}$ or $E^{\ast}_{c}$ by $E^{\ast}$, further define
$\|\cdot\|^{\ast}:E^{\ast}\rightarrow L^{0}_{+}$ by
$\|f\|^{\ast}=\bigvee\{|f(y)|~|~y\in E~and~\|y\|\leq1\}$ and
$\cdot:L^{0}(\mathcal{F},K)\times E^{\ast}\rightarrow E^{\ast}$ by
$(\xi\cdot f)(x)=\xi\cdot(f(x))$ for any $\xi\in
L^{0}(\mathcal{F},K)$, $f\in E^{\ast}$ and $x\in E$. Then it is
clear that $(E^{\ast},\|\cdot\|^{\ast})$ is an $RN$ module over $K$
with base $(\Omega,{\cal F},P)$, called the random conjugate space of $(E,\|\cdot\|)$ (see
\cite{Guo-inner}).

The following notion of a gauge function was
presented by D.Filipovi\'{c}, M.Kupper and N.Vogelpoth in \cite{FKV}
for the first time.

\no {\bf Definition 2.8 (\cite{{Guo-Relations
between},{FKV}}).}\quad {\it Let $E$ be a left module over the algebra
$L^{0}(\mathcal{F},K)$ and $A$ a subset of $E$. Then

$(1)$  $A$ is called $L^{0}-$convex if $\xi\cdot x+\eta\cdot y\in E$
for any $x$ and $y$ in $A$ and for any $\xi$ and $\eta$ in
$L^{0}_{+}$ such that $\xi+\eta=1$;

$(2)$  $A$ is called $L^{0}-$absorbent if for each $x\in E$ there
exists some $\xi\in L^{0}_{++}$ such that $x\in \xi\cdot
A:=\{\xi\cdot a~|~a\in A\}$;

$(3)$  $A$ is called $L^{0}-$balanced if $\xi\cdot x\in A$ for any
$x\in A$ and $\xi\in L^{0}(\mathcal{F},K)$ such that $|\xi|\leq 1$.}

\no {\bf Definition 2.9 (\cite{FKV}).}\quad {\it Let E be a left
module over $L^{0}(\mathcal {F}, K)$. Then the gauge function
$p_{G}:~E\rightarrow \bar{L}_{+}^{0}$ of a set $G\subset E$ is
defined by

$$p_{G}(x):=\bigwedge\{\xi\in L_{+}^{0}~|~x\in \xi\cdot G\}.$$}
\no {\bf Proposition 2.4 (\cite{FKV}).}\quad {\it Let E be a left
module over $L^{0}(\mathcal {F}, K)$. The gauge function $p_{G}$ of
an $L^{0}-$absorbent set $G\subset E$ has the following properties:

\hspace{2mm}$(\lowercase \expandafter {\romannumeral 1})$
$p_{G}(x)\leq 1$ for all $x\in G$;

\hspace{1mm}$(\lowercase \expandafter {\romannumeral 2})$
$\tilde{I}_{A}\cdot p_{G}(\tilde{I}_{A}\cdot
x)\leq\tilde{I}_{A}\cdot p_{G}(x)$ for all $x\in E$ and
$A\in\mathcal {F}$;

$(\lowercase \expandafter {\romannumeral 3})$ $\xi\cdot
p_{G}(\tilde{I}_{[\xi>0]}\cdot x)=p_{G}(\xi\cdot x)$ for all $x\in
E$ and $\xi\in L^{0}_{+}$; in particular, $\xi\cdot
p_{G}(x)=p_{G}(\xi\cdot x)$ if $\xi\in L^{0}_{++}$.}

A non-empty $L^{0}-$absorbent $L^{0}-$convex set $G\subset E$ always
contains the origin; depending on the choice of $G\subset E$, the
gauge function may be an $L^{0}-$sublinear or an $L^{0}-$seminorm.

\no {\bf Proposition 2.5 (\cite{FKV}).}\quad {\it Let E be a left
module over $L^{0}(\mathcal {F}, K)$. Then the gauge function $p_{G}$ of
an $L^{0}-$absorbent $L^{0}-$convex set $G\subset E$ satisfies:

\hspace{2.5mm}$(\lowercase \expandafter {\romannumeral 1})$
$p_{G}(x)=\bigwedge\{\xi\in L^{0}_{++}~|~x\in \xi\cdot G\}$ for all
$x\in E$;

\hspace{1.5mm}$(\lowercase \expandafter {\romannumeral 2})$
$\xi\cdot p_{G}(x)=p_{G}(\xi\cdot x)$ for all $\xi\in L_{+}^{0}$ and
$x\in E$;

$(\lowercase \expandafter {\romannumeral 3})$ $p_{G}(x+y)\leq
p_{G}(x)+p_{G}(y)$ for all $x,~y\in E$;

\hspace{0.5mm}$(\lowercase \expandafter {\romannumeral 4})$ for all
$x\in E$ there exists a sequence $\{\eta_{n}\}_{n=1}^{\infty}$ in
$L_{++}^{0}$ such that
$$\eta_{n}\searrow p_{G}(x)~~a.s.,$$
in particular, $p_{G}$ is an $L^{0}-$sublinear functional since
$0\in G$;

if $G$ is also $L^{0}-$balanced, then $p_{G}$ satisfies:

\hspace{1mm}$(\lowercase \expandafter {\romannumeral 5})$
$p_{G}(\xi\cdot x)=|\xi|\cdot p_{G}(x)$ for all $\xi\in L^{0}$ and
for all $x\in E$, namely $p_{G}$ is an $L^{0}-$seminorm.}

\no {\bf Proposition 2.6 (\cite{FKV}).}\quad {\it Let E be a left
module over $L^{0}(\mathcal {F}, K)$. Then the gauge function $p_{G}$ of
an $L^{0}-$absorbent $L^{0}-$convex set $G\subset E$ satisfies that
$p_{G}(x)\geq 1$ for all $x\in E$ with $\tilde{I}_{A}\cdot x\notin
\tilde{I}_{A}\cdot G$ for all $A\in\mathcal {F}$ with $P(A)>0$.}


\section{The geometric form of Hahn-Banach extension theorem for $L^{0}-$linear functions}

\no{\bf Theorem 3.1}\quad Proposition 1.1 is equivalent to Proposition 1.2.

\no {\bf Proof.}\quad Let $p_{G}$ be the gauge function of $G$,
namely $$p_{G}(x)=\bigwedge\{\xi\in L_{++}^{0}~|~\xi\cdot x\in
G\},~\forall x\in E.$$ Since $G$ is an $L^{0}-$convex and $L^{0}-$absorbent subset
of $E$, $p_{G}$ is an $L^{0}-$sublinear functional on $E$ by
Proposition 2.5, and since $g(y)\leq 1$ for any $y\in M\cap G$, then
for any $x\in M$ and $\lambda\in L_{++}^{0}$ we can obtain
$g(x)\leq \lambda$ when $x\in \lambda\cdot G$, namely $g(x)\leq
p_{G}(x)$. From Proposition 3.1, there exists an $L^{0}-$linear
functional $f:E\rightarrow L^{0}({\cal F},R)$ such that $f$ extends
$g$ and $f(x)\leq p_{G}(x),\forall x\in E$. Therefore, we have that
$$f(x)\leq 1,\forall x\in G.$$

Conversely, let $G=\{x\in E~|~p(x)\leq 1\}$, then it is clear that $g,~M$ and $G$
satisfy Proposition 1.2, hence there exists an $L^{0}-$linear
functional $f:E\rightarrow L^{0}({\cal F},R)$ such that $f$ extends
$g$ and
$$f(x)\leq 1,\forall x\in G$$ by Proposition 1.2, so that we can have that
$f(x)\leq p(x),\forall x\in E$. $\square$

If $K=C$, we have the following geometric form:

\no{\bf Theorem 3.2.}\quad Let $E$ be a left module over the algebra
$L^{0}(\mathcal{F},C)$, $M$ an $L^{0}(\mathcal{F},C)-$submodule in
$E$ and $G$ an $L^{0}-$convex and $L^{0}-$absorbent subset of $E$.
If $g:M\rightarrow L^{0}({\cal F},C)$ is an $L^{0}-$linear
functional and $(Reg)(y)\leq 1$ for any $y\in M\cap G$, then there
exists an $L^{0}-$linear functional $f:E\rightarrow L^{0}({\cal
F},C)$ such that $f$ extends $g$ and $(Ref)(x)\leq 1,\forall x\in
G$.

Now, we give a new proof of Proposition 1.3 by the geometric form of
Hahn-Banach extension theorem for $L^{0}-$functions. In fact, we
need only to prove the following basic strict separation theorem for the case of
$RN$ modules, namely Proposition 3.1 below, since by Proposition 3.1 one can easily complete the remaining part of the proof of Proposition 1.3, see
\cite{GXC} for details.

Let $(E,\|\cdot\|)$ be an $RN$ module, $G$ a
$\mathcal{T}_{\varepsilon,\lambda}-$closed $L^{0}-$convex subset of
$E$ and $x_{0}\in E\setminus G$. Then $\{\|x_{0}-g\|~|~ g\in G\}$ is
a dually directed subset in $L^{0}_{+}$ and one can obtain the
following claim: there exists  $A\in \mathcal {F}$ with $P(A)>0$
such that $\bigwedge\{\|x_{0}-g\|~|~ g\in G\}>0$ on $A$ from Lemma
3.8 in \cite{Guo-Relations between}.

\no{\bf Proposition 3.1 (\cite{GXC}).}\quad {\it Let $(E,
\|\cdot\|)$ be an RN module over $R$ with base $(\Omega,\mathcal
{F},P)$, $G$ a $\mathcal {T}_{\varepsilon,\lambda}-$closed and
$L^{0}-$convex subset of $E$, $x_0\in E\setminus G$, and $
\xi=\bigwedge\{\|x_{0}-h\|~|~h\in G\}$.
Then there exists a continuous module homomorphism $f$ from
$(E,\mathcal {T}_{\epsilon,\lambda})$ to $(L^{0}(\mathcal
{F},R),\mathcal {T}_{\epsilon,\lambda})$ such that
$$f(x_0)>\bigvee\{|f(y)|~|~ y\in G\}$$ and $$f(x_0)>\bigvee\{|f(y)|~|~ y\in
G\}~ on~[\xi>0].$$ }

\no {\bf Proof.} Without loss of generality, we can assume
$\theta\in G$ (otherwise, by a translation). It is easy to see that
$\tilde{I}_{B}\cdot x_{0}\not\in G$ for any $B\in \mathcal {F}$ with
$B\subset A$ and $P(B)>0$. In fact, assume that
$\tilde{I}_{B}\cdot x_{0}\in G$, then $\xi=0$ on $B$, which is
contradict to $\tilde{A}=[\xi>0]$ and $B\subset A$. Let $M=\{x\in
E~|~ \|x\|\leq \frac{1}{3}\tilde{I}_{A}\cdot \xi~on~A\}$, then it is
clear that $M$ is $L^{0}-$convex and $L^{0}-$absorbent. Further, let
$G+M=\{h+x~|~ h\in G, x\in M\}$, then $G+M$ is also an
$L^{0}-$convex and $L^{0}-$absorbent subset of $E$. Since $\theta\in
G+M$ and $G+M$ is an $L^{0}-$convex, we have that
$\tilde{I}_{F}\cdot (G+M)\subset G+M$ for every subset
$F\in\mathcal {F}$ with $P(F)>0$.

Let $p_{G+M}$ be the gauge function of $G+M$, then $p_{G+M}$ is an
$L^{0}-$sublinear functional on $E$ by Proposition 2.5. It is easy
to see that $p_{G+M}(x)=0$ on $A^{c}$ for any $x\in E$. Now we prove
that $p_{G+M}(x_{0})> 1$ on $A$. In fact, for any $z=z_{G}+z_{M}$,
where $z_{G}\in G$ and $z_{M}\in M$, since
$\|(x_{0}-z)\|\geq\|x_{0}-z_{G}\|-\|z_{M}\|$, we can obtain that
$\|(x_{0}-z)\|\geq\frac{2}{3}\xi>0$ on $A$ from
$\|x_{0}-z_{G}\|\geq\xi$ and $\|z_{M}\|\leq\frac{1}{3}\cdot\xi$ on
$A$. Thus $x_{0}\not\in G+M$ and $p_{G+M}(x_{0})\neq 0$ by
Definition 2.8. From Proposition 2.5, there exists a sequence
$\{\eta_{n}~|~n\in N\}\subset L_{++}^{0}$ such that $x_{0}\in
\eta_{n}\cdot (G+M)$ and $\eta_{n}\searrow p_{G+M}(x_{0})$.
According to $x_{0}\in E\setminus (G+M)$ and
$\bigwedge\{\|x_{0}-h\|~|~ h\in G+M\}>\frac{1}{3}\xi$ on $A$, we
have that $\eta_{n}>1$ on $A$ for any $n\in N$ and hence
$p_{G+M}(x_{0})\geq 1$ on $A$. Let
$\tilde{D}=[p_{G+M}(x_{0})=1]\cap\tilde{A}$, then we will prove that
$P(D)>0$ is impossible: if $P(D>0)$, it is clear that
$\eta_{n}\searrow \tilde{I}_{D}$ and $Q(\eta_{n})\nearrow
\tilde{I}_{D}$ on $D$; since $x_{0}\in \eta_{n}\cdot(G+M)$ and
$\eta_{n}>1$ on $D$, we have $Q(\eta_{n})\cdot x_{0}\in
(Q(\eta_{n})\cdot\eta_{n})\cdot (G+M)\subset G+M$ and
$\|x_{0}-Q(\eta_{n})\cdot x_{0}\|=\|(1-Q(\eta_{n}))\cdot
x_{0}\|\searrow 0$ on $D$, which contradicts to the fact that
$\bigwedge\{\|x_{0}-h\|~|~ h\in G+M\}>\frac{1}{3}\xi$ on $A$. Hence
$P(D)=0$ and $p_{G+M}(x_{0})>1$ on $A$.

Now we prove that the gauge function $p_{G+M}$ of $G+M$ is
continuous under the $(\epsilon,\lambda)-$topology. Let $x$ be an
arbitrary element of $E$, $\tilde{H}_{x}=[\|x\|\neq 0]$ and
$\overline{t}=\frac{1}{3}\tilde{I}_{A\cap H_{x}}\cdot \xi\cdot
\|x\|^{-1}$, then $\overline{t}\cdot x\in M\subset G+M$. Thus we
have
$$p_{G+M}(x)=0 $$
on $H_{x}^{c}$ and
$$p_{G+M}(x)\leq 3\tilde{I}_{A\cap H_{x}}\cdot Q(\xi)\cdot \|x\|$$
on $H_{x}$.

\no Therefore we can obtain that
$$p_{G+M}(x)\leq (3Q(\xi)+1)\cdot \|x\| $$
and $p_{G+M}$ is continuous under the $(\epsilon,\lambda)-$topology.

Let $U=\{k\cdot x_{0}~|~k\in L^{0}(\mathcal{F},R)\}$, then $U$ is an
$L^{0}-$submodule of $E$. Define an $L^{0}-$linear function
$g:U\rightarrow L^{0}(\mathcal {F}, R)$, where
$g(x_{0})=\tilde{I}_{A}\cdot p_{G+M}(x_{0})$. Then we have that
$g(y)\leq 1,~\forall
y\in U\cap (G+M)$. By Proposition 1.2, there exists an
$L^{0}-$linear function $\bar{g}:E\rightarrow L^{0}({\cal F},R)$
such that $\bar{g}$ extends $g$ and $\bar{g}(x)\leq 1,\forall y\in
G+M$. Hence, $\bar{g}(x_{0})>1$ on $A$ and
$$\bar{g}(x_{0})>\bigvee\{|f(y)|~|~ y\in G\} ~on~A.$$ Let $f=\tilde{I}_{A}\cdot \bar{g}$, since $f\leq
p_{G+M}$ from Remark 3.1, it is easy to see that $f\in
E^{\ast}_{\varepsilon,\lambda}$ and $f(x)=0$ on $A^{c}$. Therefore,
it is clear that
$$f(x_{0})>\bigvee\{|f(y)|~|~ y\in G\}$$ and
$$f(x_{0})>\bigvee\{|f(y)|~|~ y\in G\} ~on~A.~~~\square$$


\section{The Goldstine-Weston theorem in $RN$ modules}

Before giving the proofs of Theorem 1.1 and 1.2, we first present
some necessary definitions and lemmas.

\no {\bf Definition 4.1 (\cite{Guo-Relations between,
Guo-progress}).}\quad {\it Let $E$ be a left module over the algebra
$L^{0}({\cal F},K)$. Such a formal sum $\sum_{n\geq
1}\tilde{I}_{A_n}x_n$ for some countable partition $\{A_n,n\in N\}$
of $\Omega$ to ${\cal F}$ and some sequence $\{x_{n}~|~n\in N\}$ in
$E$, is called a countable concatenation of $\{x_{n}~|~n\in N\}$
with respect to $\{A_n,n\in N\}$. Furthermore a countable
concatenation $\sum_{n\geq 1}\tilde{I}_{A_n}x_n$ is well defined or
$\sum_{n\geq 1}\tilde{I}_{A_n}x_n\in E$ if there is $x\in E$ such
that $\tilde{I}_{A_{n}}x=\tilde{I}_{A_{n}}x_{n},~\forall n\in N$. A
subset $G$ of $E$ is called having the countable concatenation
property if every countable concatenation $\sum_{n\geq
1}\tilde{I}_{A_n}x_n$ with $x_{n}\in G$ for each $n\in N$ still
belongs to $G$, namely $\sum_{n\geq 1}\tilde{I}_{A_n}x_n$ is well
defined and there exists $x\in G$such that $x=\sum_{n\geq
1}\tilde{I}_{A_n}x_n$. }

\no {\bf Lemma 4.1 (\cite{Guo-Relations between}).} \quad {\it Let
$(E,\mathcal{P})$ be a random locally convex module over $K$ with
base $(\Omega,\mathcal{F},P)$, $G\subset E$ a subset having the
countable concatenation property. Then
$\bar{G}_{\varepsilon,\lambda}=\bar{G}_{c}$. }

Now, let us recall the random weak topology and the random weak star topology.

\no  {\bf  Definition 4.2 (\cite{G-C,Guo-progress}).}\quad {\it Let
$(E,\|\cdot\|)$ be an $RN$ module over $K$ with base
$(\Omega,\mathcal{F},P)$, $(E^{\ast},\|\cdot\|^{\ast})$ the random
conjugate space of $E$. For any $f\in E^{\ast}$, define $\|\cdot\|_{f}:E\rightarrow L^{0}_{+}$ by $\|x\|_{f}=|f(x)|,~\forall x\in E$, and denote $\{\|\cdot\|_{f}~|~f\in E^{\ast}\}$ by $\sigma(E,E^{\ast})$, it is clear that $(E,\sigma(E,E^{\ast}))$ is a random locally convex module over $K$ with base
$(\Omega,\mathcal{F},P)$. Then the $(\varepsilon,\lambda)-$topology $\sigma_{\varepsilon,\lambda}(E,E^{\ast})$ and the locally $L^{0}-$convex topology $\sigma_{c}(E,E^{\ast})$ on $E$ induced by $\sigma(E,E^{\ast})$ are called random weak $(\varepsilon,\lambda)-$topology and random weak locally $L^{0}-$convex topology on $E$, respectively. }

\no  {\bf  Remark 4.1.}\quad Similarly, we can define the random
weak star $(\varepsilon,\lambda)-$topology
$\sigma_{\varepsilon,\lambda}(E^{\ast},E)$ and the random weak star locally
$L^{0}-$convex topology $\sigma_{c}(E^{\ast},E)$ on $E^{\ast}$, respectively.

\no {\bf Lemma 4.2 (\cite{G-C}).} \quad {\it Let $(E,\|\cdot\|)$ be
an $RN$ module over $K$ with base $(\Omega,\mathcal{F},P)$, then $(E^{\ast},\sigma_{c}(E^{\ast},E))^{\ast}=E$. Furthermore, if
$E$ has the countable concatenation property, then
 $(E^{\ast},\sigma_{\varepsilon,\lambda}(E^{\ast},E))^{\ast}
=E$.}

If $(B,\|\cdot\|)$ is a normed space and
$(B^{\prime},\|\cdot\|^{\prime})$ is the classical conjugate space of $B$, we
have that $B^{\prime}(1)=\{f\in B^{\prime}~|~\|f\|^{\prime}\leq 1\}$
is compact under the weak star topology by the well known Banach-Alaoglu
theorem. Hence, $B^{\prime}(1)$ is closed with respect to the weak
star topology. Let $(E,\|\cdot\|)$ be an $RN$ module over $K$ with
base $(\Omega,\mathcal{F},P)$ and $(E^{\ast},\|\cdot\|^{\ast})$ the
random conjugate space. In \cite{Guo-BBKS}, Guo proved that
$E^{\ast}(1)=\{f\in E^{\ast}~|~\|f\|^{\ast}\leq 1\}$ is not compact
under $\sigma_{c}(E^{\ast},E)$ unless $(\Omega,\mathcal{F},P)$ is
essentially purely $P-$atomic. But
Lemma 4.3 below indicates that $E^{\ast}(1)$ is still closed with respect to both
$\sigma_{\varepsilon,\lambda}(E^{\ast},E)$ and
$\sigma_{c}(E^{\ast},E)$.

\no {\bf Lemma 4.3.} \quad {\it Let $(E,\|\cdot\|)$ be an $RN$
module over $K$ with base $(\Omega,\mathcal{F},P)$ such that $E$ has
the countable concatenation property. Then $E^{\ast}(1)=\{f\in
E^{\ast}~|~\|f\|^{\ast}\leq 1\}$ is closed with respect to both
$\sigma_{\varepsilon,\lambda}(E^{\ast},E)$ and
$\sigma_{c}(E^{\ast},E)$. }

\no {\bf Proof.}\quad Since it is clear that $E^{\ast}(1)$ has the
countable concatenation property, we need only to check that
$E^{\ast}(1)$ is closed with respect to $\sigma_{c}(E^{\ast},E)$.
For any $f\in E^{\ast}\setminus E^{\ast}(1)$, there exists
$A\in\mathcal{F}$ such that $P(A)>0$ and $\|f\|^{\ast}>1$ on $A$. From
$\|f\|^{\ast}=\bigvee\{|f(x)|~|~\|x\|\leq1 \}$, there are $x_{f}\in
E$, $\|x_{f}\|\leq 1$ and $B\in\mathcal{F}$, $B\subset A$ with
$P(B)>0$ such that $|f(x_{f})|>1$ on $B$. Let
$$\varepsilon=\tilde{I}_{B^{c}}+\frac{|f(x_{f})|-1}{2}\cdot\tilde{I}_{B}$$
and $$B(x_{f},\varepsilon)=\{g\in E^{\ast}~|~|g(x_{f})|\leq
\varepsilon\},$$ then $B(x_{f},\varepsilon)$ is a neighborhood of
$\theta$ in $E^{\ast}$ with respect to $\sigma_{c}(E^{\ast},E)$ and,
for any $h\in B(x_{f},\varepsilon)$ it is easy to see that
$$|(f+h)(x_{f})|\geq|f(x_{f})|-|h(x_{f})|$$ and
$$(f+h)(x_{f})|\geq\frac{|f(x_{f})|+1}{2}>1$$ on $B$. Hence,
$f+h\not\in E^{\ast}(1)$, namely $f+B(x_{f},\varepsilon)\subset
E^{\ast}\setminus E^{\ast}(1)$. Consequently, $E^{\ast}(1)$ is
closed with respect to $\sigma_{c}(E^{\ast},E)$. $\square$

\no {\bf Lemma 4.4} \quad {\it Let $(E,\|\cdot\|)$ be an $RN$ module
over $K$ with base $(\Omega,\mathcal{F},P)$ such that $E$ has the
countable concatenation property, $J$, $E(1)$, $J(E(1))$ and
$\overline{J(E(1))}_{\varepsilon,\lambda}^{w^{\ast}}$ the same as in
Theorem 1.1. Then
$\overline{J(E(1))}_{\varepsilon,\lambda}^{w^{\ast}}=E^{\ast\ast}(1)$.
}

\no {\bf Proof.} \quad Since $E^{\ast\ast}(1)=\{\varphi\in
E^{\ast\ast}(1)~|~\|\varphi\|^{\ast\ast}\leq 1\}$ is closed with
respect to $\sigma_{\varepsilon,\lambda}(E^{\ast\ast},E^{\ast})$
from Lemma 4.3, it follows that
$\overline{J(E(1))}_{\varepsilon,\lambda}^{w^{\ast}}\subset
E^{\ast\ast}(1)$.

Now, we prove that $E^{\ast\ast}(1)\subset
\overline{J(E(1))}_{\varepsilon,\lambda}^{w^{\ast}}$. We only need
to prove that for any $\psi\in E^{\ast\ast}\setminus
\overline{J(E(1))}_{\varepsilon,\lambda}^{w^{\ast}}$ there is
$A_{\psi}\in\mathcal{F}$ such that $P(A_{\psi})>0$ and
$\|\psi\|^{\ast\ast}>1$ on $A_{\psi}$. Since $E$ has the countable
concatenation property, we have that
$(E^{\ast\ast},\sigma_{\varepsilon,\lambda}(E^{\ast\ast},E^{\ast}))^{\ast}
=E^{\ast}$ by Lemma 4.2 and that there exists $\bar{f}\in E^{\ast}$ such
that
$$(\mbox{Re}\bar{f})(\psi)>\bigvee\{(\mbox{Re}\bar{f})(g)~|~ g\in \overline{J(E(1))}_{\varepsilon,\lambda}^{w^{\ast}}\}$$ and
$$(\mbox{Re}\bar{f})(\psi)>\bigvee\{(\mbox{Re}\bar{f})(g)~|~ g\in \overline{J(E(1))}_{\varepsilon,\lambda}^{w^{\ast}}\}~on~[\xi>0]$$ by
Proposition 1.3, where $\xi$ is the same as in Proposition 1.3. For
any $y\in \overline{J(E(1))}_{\varepsilon,\lambda}^{w^{\ast}}$, it is easy to see that $||\bar{f}(y)|\cdot Q(\bar{f}(y))|\leq 1$ on $\Omega$,
$(|\bar{f}(y)|\cdot Q(\bar{f}(y)))\cdot y\in
\overline{J(E(1))}_{\varepsilon,\lambda}^{w^{\ast}}$ and
$\bar{f}((|\bar{f}(x)|\cdot Q(\bar{f}(y)))\cdot y)=|\bar{f}(y)|$.
Hence, we have that $$\bigvee\{(\mbox{Re}\bar{f})(g)~|~
g\in
\overline{J(E(1))}_{\varepsilon,\lambda}^{w^{\ast}}\}=\bigvee\{(|\bar{f}(g)|~|~
g\in \overline{J(E(1))}_{\varepsilon,\lambda}^{w^{\ast}}\}.$$ Let
$f=Q(\|\bar{f}\|^{\ast})\cdot \bar{f}$ and $A_{\psi}=[\xi>0]$, then
we have that $\|f\|^{\ast}=1$ on $A_{\psi}$ and
$$(\mbox{Re}f)(\psi)>\bigvee\{(\mbox{Re}f)(g)~|~ g\in
\overline{J(E(1))}_{\varepsilon,\lambda}^{w^{\ast}}\}$$ and
$$(\mbox{Re}f)(\psi)>\bigvee\{(\mbox{Re}f)(g)~|~ g\in \overline{J(E(1))}_{\varepsilon,\lambda}^{w^{\ast}}\}~on~A_{\psi}.$$
Consequently, we can obtain $$\|\psi\|^{\ast\ast}\geq|f(\psi)|\geq
(\mbox{Re}f)(\psi)>\bigvee\{(\mbox{Re}f)(g)~|~ g\in
\overline{J(E(1))}_{\varepsilon,\lambda}^{w^{\ast}}\}$$
~~~~~~~~~~~~~~~~~=$\bigvee\{|f(g)|~|~ g\in
\overline{J(E(1))}_{\varepsilon,\lambda}^{w^{\ast}}\}\geq\bigvee\{|f(y)|~|~
y\in E(1)\}=\|f\|^{\ast}$ on $A_{\psi}$,

\no namely $\|\psi\|^{\ast\ast}>1$ on $A_{\psi}$. $\square$

\no  {\bf  Definition 4.4 (\cite{Guo-Relations between}).} \quad
{\it Let $E$ be a left
module over the algebra $L^{0}(\mathcal{F},K)$ and $G$ a subset of $E$. The set of countable concatenations
$\Sigma_{n\geq1}\tilde{I}_{A_{n}}x_{n}$ with $x_{n}\in G$ for
each $n\in N$ is called the countable concatenation hull of $G$, denoted by $H_{cc}(G)$.  }

\no {\bf Proof of Theorem 1.1.}\quad Denote $H_{cc}(E)$ by $E_{cc}$ and
define $\|\cdot\|_{cc}:E_{cc}\rightarrow L^{0}_{+}$ by
$\|x\|_{cc}=\sum_{n\leq1}\tilde{I}_{A_{n}}\cdot \|x_{n}\|$ for any
$x=\sum_{n\leq1}\tilde{I}_{A_{n}}\cdot x_{n}$ in $E_{cc}$, where
$\{A_{n}~|~n\in N\}$ is a countable partition of $\Omega$ to
$\mathcal{F}$ and $x_{n}\in E$ for any $n\in N$. It is easy to see
that $E^{\ast\ast}_{cc}=E^{\ast\ast}$. By Theorem 4.1, we can
obtain that
$\overline{J(E_{cc}(1))}_{\varepsilon,\lambda}^{w^{\ast}}=E^{\ast\ast}(1)$.
Since $J(E(1))$ is dense in $J(E_{cc}(1))$ with respect to the
$(\varepsilon,\lambda)-$topology which is induced by
$\|\cdot\|_{cc}$ and stronger than
$\sigma_{\varepsilon,\lambda}(E^{\ast\ast},E^{\ast})$,
our desired result follows from the fact that the $(\varepsilon,\lambda)-$topology is stronger than
$\sigma_{\varepsilon,\lambda}(E^{\ast\ast},E^{\ast})$. $\square$

\no {\bf Proof of Theorem 1.2.}\quad It follows immediately from
Theorem 1.1 and Lemma 4.1. $\square$

The following example shows that $J(E(1))$ may be not dense in
$E^{\ast\ast}(1)$ under $\sigma_{c}(E^{\ast\ast},E^{\ast})$ if $E$
has not the countable concatenation property.

\no {\bf Example 4.1.} \quad {\it Let $\Omega=\{1,2,3,\cdots \}$,
$\mathcal {F}=2^{\Omega}$, $\bar{P}:\mathcal {F}\rightarrow R$ such
that $\bar{P}(\Lambda)=$ the number of points in $\Lambda$ if
$\Lambda$ is any finite subset in $\Omega$ and
$\bar{P}(\Lambda)=\infty$ otherwise and $P:\mathcal {F}\rightarrow
[0,1]$ such that
$P(\Lambda)=\sum_{n=1}^{\infty}\frac{1}{2^{n}}\frac{\bar{P}(\Lambda\cap
\{n\})}{\bar{P}(\{n\})}$ for each subset $\Lambda$ of $\Omega$, then $(\Omega,\mathcal {F},P)$ is a
probability space. Let
$(E,\|\cdot\|)=(L^{0}(\mathcal{F},K),|\cdot|)$ and $F=\{\varphi\in
E~|~$there is a positive integer $n_{\varphi}$ such that
$\varphi(k)=0,~\forall k\geq n_{\varphi}\}$, then it is clear that $F$ is an $L^{0}-$submodule of $E$ and $(F,|\cdot|)$ is
an $RN$ module over $K$ with the base $(\Omega,\mathcal {F},P)$. Let $F(1)=\{x\in F~|~|x|\leq 1\}$, then
it is easy to check that $F(1)$ is a closed subset in both $(E,|\cdot|)$ and
$(F,|\cdot|)$ under $\mathcal {T}_{c}$ induced by $|\cdot|$. Hence,
we have that $F(1)$ is not dense in $E(1)$ under $\mathcal {T}_{c}$.
Furthermore, it is clear that
$(F^{\ast},\|\cdot\|^{\ast})=(F^{\ast\ast},\|\cdot\|^{\ast\ast})=(E,|\cdot|)$
and $\sigma_{c}(F^{\ast\ast},F^{\ast})$ is also the locally $L^{0}-$
convex topology induced by $|\cdot|$. Consequently, $J(F(1))$ is not
dense in $F^{\ast\ast}(1)$ under
$\sigma_{c}(F^{\ast\ast},F^{\ast})$. }

\no {\bf Acknowledgements} ~~~~The authors would like to thank our
supervisor Professor Guo TieXin for some helpful and critical suggestions which considerably improve the readability of this paper.

\footnotesize
\parindent=6mm
\makeatletter
\renewcommand\section{\@startsection {section}{1}{\z@}%
{-3.5ex \@plus -1ex \@minus -.2ex}%
{2.3ex \@plus.2ex}{\bfseries\large}}%

\makeatother

\end{document}